\documentclass[11pt,reqno,a4paper]{article}
\usepackage{amsmath,amsthm,amsfonts,amssymb}

\newtheorem{tm}{Theorem}[section]

\newtheorem{prop}[tm]{Proposition}

 \theoremstyle{definition}

\newcommand{\beqa}{\begin{eqnarray*}}
\newcommand{\eeqa}{\end{eqnarray*}}

 \def\cF{\mathcal{F}}              

 \def\cH{\mathcal{H}}
 \def\cB{\mathcal{B}}
 \def\cE{\mathcal{E}}
 \def\cG{\mathcal{G}}

 \def\cA{\mathcal{A}}

\def\L{\left(}
\def\R{\right)}

\def\<{\left<}
\def\>{\right>}

\def\inv{^{-1}}

\def\mv1{M_v^1}

\begin{document}
\title{On some stability results of frame atomic decompositions
\thanks{The author acknowledges the partial support provided through the
FP6 Intra-European Individual Marie Curie Fellowship Programme, project FTFDORF-501018.}}
\author{Massimo~Fornasier}
\maketitle
\begin{center}
{\it In honor of Prof. Laura Gori for her 70$^\text{\it th}$ birthday}
\end{center}
\begin{abstract} This paper is concerned with the implications of sufficient conditions ensuring that a perturbation of a frame is again a frame. We emphasize how stability of frames is  fundamental for numerical applications and we discuss in particular the connection between stability conditions and localization principles for frame atomic decompositions in Banach spaces.

\end{abstract}

\noindent
{\bf AMS subject classification:} 42C15, 46B99, 46H99, 65D99, 65J05, 65T60

\noindent
{\bf Key Words:} Banach spaces, Banach algebras, localization of frames, non-orthogonal expansions, numerical stability, perturbations.
\section{Introduction}

Stable redundant non-orthogonal expansions in Hilbert spaces has been introduced  by Duffin and Schaeffer \cite{DS} under the name of {\it frames}. 
Besides traditional and relevant applications of {\it wavelet and Gabor frames} \cite{DGM,D1,D2,G,FS,FS1} in signal processing, image processing, data compression,   pattern matching, sampling theory, communication and data transmission, recently the use of frames also in numerical analysis for the solution of operator equations and PDE is investigated \cite{S,DFR1}.
Therefore, not only the characterization  by frames of functions in $L^2(\mathbb{R}^d)$ is relevant but also that of (smoothness) Banach  function spaces is crucial to have a correct formulation of effective and stable numerical schemes. The concept of {\it frame atomic decomposition} as an extension of frames in {\it coorbit spaces} has been introduced in \cite{FG3,G0} and it will be an important concept in this paper.
The problem when one is dealing with numerical applications is that one cannot use exact representations of functions, but only approximations can be available.
Therefore any frame that one wants concretely use in applications will be affected by a perturbation.
At this point it is not ensured that the perturbed system is a frame anymore.
For example, recent papers \cite{CS,DFR1} proposed efficient methods to compute canonical dual frames. Such systems are fundamental for the computation of the frame expansion coefficients by scalar product.
In fact the proposed algorithms compute a nice and accurate approximation of each individual element of the canonical dual frame, producing a global system that can be interpreted as a perturbation of the true and original canonical dual.
On the other hand, no proof is given yet whether such new system is again a frame and how much the perturbation produced by the numerics can affect global expansions. 
\\

In this paper we want to discuss several stability results of frames under different kinds of perturbations, combining some of the well-known results due to Casazza, Christensen, and Heil \cite{C2,cc1,CH,FS} with some new insights in the frame theory. In particular, the emphasis here will be on {\it localization of frame principles} introduced by Gr\"ochenig in \cite{G1} and on their generalizations due to Gr\"ochenig and the author \cite{FG}. We will show that a suitable perturbation of an intrinsically localized frame can produce again an intrinsically localized frame atomic decomposition  for a class of associated Banach spaces. Applications of such result will be presented in subsequent contributions.

The paper is organized as follows: Section 2 recalls the concepts of frames and some relevant instances of the localization of frame theory and Section 3 collects the perturbation and stability results.

\section{Frames and Schur localization}

In this section we recall the concept of frames, how they  can be used to define certain associated
Banach spaces, and how to obtain  stable decompositions in these  Banach spaces.


A subset  $\cG=\{g_n\}_{n \in {\mathbb{Z}}}$ of a separable Hilbert
 space   $\mathcal{H}$ is called \emph{frame} for $\mathcal{H}$ if
\begin{equation}
 A \|f\|^2 \leq \sum_{n \in {\mathbb{Z}}} |\langle f, g_n \rangle|^2 \leq B \|f\|^2, \quad \forall f\in \mathcal{H},
\end{equation}
for some  constants $0<A \leq B < \infty$.

Equivalently, we could define a frame by the requirement that  the frame
operator  $S=D  C$ is a boundedly invertible (positive and
self-adjoint) on $\cH $, where $C=C_{\cG }$ defined by
$Cf = (\langle f, g_n \rangle)_n $ is the 
corresponding bounded \emph{analysis   operator} from $ \mathcal{H}
\rightarrow \ell^2({\mathbb{Z}})$, and $D=
D_{\cG } = C^*, D\mathbf{c} = \sum_n c_n g_n, $  is the  bounded \emph{synthesis operator} from $
\ell^2({\mathbb{Z}}) \rightarrow \mathcal{H}$. The set  $\tilde{\cG}=S^{-1} \cG$ is again a
frame for $\mathcal{H}$ and it is called the  \emph{canonical dual frame}
playing an important role in the reconstruction of $f\in \cH $ from the
frame coefficients, because we have
\begin{equation}
f =S S^{-1} f = \sum_n \langle f, S^{-1} g_n \rangle g_n=S^{-1} S f = \sum_n \langle f, g_n \rangle S^{-1} g_n.
\end{equation}
Since in general  a frame is  overcomplete,  the coefficients in this
expansion are  in general not unique (unless $\cG $ is a Riesz basis, we have
$\text{ker}(D) \neq \{0\}$) and  there exist many possible dual frames  $\{
\tilde g_n \}_{n \in {\mathbb{Z}}}$ in $\mathcal{H}$ such that
$$
        f = \sum_{n \in {\mathbb{Z}}} \langle f, \tilde g_n \rangle g_n
$$
with the norm equivalence $\|f\|_{\cH } \asymp \|\langle f, \tilde g_n \rangle
\|_2$.
We refer the reader, for example, to the book \cite{C} for a more complete literature, details, and examples.

The concept of frame can be extended to Banach spaces as follows: A \emph{frame atomic decomposition} for a separable Banach space $\cB$ is a sequence $\cG =\{g_n\}_{n \in \mathbb{Z}}$ in $\cB$ with an associated sequence space $\cB_d$ such that the following properties hold.
\begin{itemize}
\item[(a)] There exists a \emph{coefficient operator} $C$ defined by $C f=\L \langle f,\tilde g_n\rangle_{n \in {\mathbb{Z}}} \R$ bounded from $\cB$ into $\cB_d$, where $\tilde \cG =\{\tilde g_n\}_{n \in {\mathbb{Z}}}$ is in $\cB'$;
\item[(b)] norm equivalence: for all $f \in \cB$
$$
        A_{\cB} \|f\|_{\cB} \leq \| \langle f,\tilde g_n\rangle_{n \in {\mathbb{Z}}}\|_{\cB_d} \leq  B_{\cB}\|f\|_{\cB}, \quad   A_{\cB}, B_{\cB}>0;
$$
\item[(c)] the following series expansion converge unconditionally
$$
        f= \sum_{n \in \mathbb{Z}} \langle f,\tilde g_n\rangle g_n, \quad \text{for all } f\in \cB.
$$
\end{itemize}
We want to illustrate in the following that there exists a natural choice of Banach spaces $\cB$ such that ``suitable'' frames $\cG$ for $\cH$ extend to frame atomic decompositions for $\cB$ with coefficient map $C=C_{\tilde \cG}$ where $\tilde \cG$ denotes the canonical dual frame.
For ``suitable'' we mean that the frame should have an additional {\it localization structure}. \\

The theory of localization of frames has been introduced and developed by Gr\"ochenig {\it et al.} \cite{G1,G2,CG,FG} in order to illustrate general principles for the frame characterization of Banach spaces.
Recently has been also recognized that localized frames are ``good'' not only for theoretical purposes but also for numerical applications in signal and image processing and in PDE \cite{FS,FS1,CS,S,DFR1}.

In this paper we work only with a particular type of localization and we refer to \cite{FG} for a more general theory and results.
In particular we shall work with the \emph{Schur algebra} $\cA_s^1$ \cite{GL} which is defined as   the class of matrices $A= (a_{kl}), k,l \in \mathbb{Z}  $, such
that
$$
\|A\|_{\cA_s^1}:=\max\{\sup_{k \in \mathbb{Z}} \sum_{l \in \mathbb{Z}} |a_{kl}|v_s(k-l), \sup_{l \in \mathbb{Z}} \sum_{k \in \mathbb{Z}} |a_{kl}| v_s(k-l)\} < \infty,
$$
where $v_s(x)=(1+|x|)^s$ for $s \geq 0$. For $s=0$ we denote $\cA^1_s$ with $\cA^1$.
The Schur algebra endowed with the norm $\|\cdot\|_{\cA_s^1}$ is a Banach $*$-algebra, where the involution is the transpose-conjugate operator and the following properties hold

\begin{itemize}
\item[(A0)] $\cA^1_s \subseteq \cA^1 \subseteq \cB (\ell ^p({\mathbb{Z}} ) )$, i.e., each $A \in
  \cA^1_s $ defines a bounded operator on $\ell ^p ({\mathbb{Z}} )$ for $1 \leq p \leq \infty$;
\item[(A1)] If $s>0$ and $A \in \cA^1_s$ is invertible on $\ell ^2({\mathbb{Z}} )$, then
  $A\inv \in \cA^1_s $ as well. In the language of Banach algebras, $\cA^1_s $
  is called inverse-closed in $\cB (\ell ^2({\mathbb{Z}} ) ) $. It is not known whether such property holds even for $s=0$;
\item[(A2)] $\cA^1_s $ is solid:  i.e., if $A\in \cA^1_s $ and  $|b_{kl}| \leq
  |a_{kl}|$ for all $k,l\in {\mathbb{Z}} $,  then  $B\in \cA^1_s$ as well.
\end{itemize}
 
We refer to  \cite{GL} where a characterization of a large class of algebras with properties (A0-2) is presented. 
We also observe here that a major part of the results presented in this paper can be generalized considering those algebras instead of  $\cA_s^1$.

Given two frames $\cG=\{g_n\}_{n \in {\mathbb{Z}}}$ and $\mathcal{F}=\{f_x\}_{x \in {\mathbb{Z}}}$ for the Hilbert space $\mathcal{H}$, the (cross-) Gramian matrix $A=A(\cG,\cF)$ of $\cG $ with respect to $\cF $ is the $\mathbb{Z} \times {\mathbb{Z}} $-matrix with entries
$$
a_{nx} = \langle g_n, f_x \rangle .
$$
A frame $\cG $ for $\cH $ is called \textit{$\cA^1_s$-localized} with respect to
another frame $\cF $ if $A(\cG,\cF) \in \cA^1_s$. In this case we write $\cG \sim
_{\cA^1_s } \cF $.  If $\cG \sim_{\cA^1_s } \cG$, then $\cG$ is called \textit{$\cA^1_s $-self-localized} or
\textit{intrinsically $\cA^1_s$-localized}. By exploiting property (A1) one can prove \cite{FG} the following 

\begin{tm}
\label{self}
For $s>0$, any  $\cA^1_s$-self-localized frame $\mathcal{G}$ has always $\cA^1_s$-self-localized  canonical dual.
\end{tm}

Let $(\cG,\tilde \cG)$ be a pair of dual  $\cA^1_s$-self-localized frames for $\mathcal{H}$. 
Assume $1 \leq p\leq 2$. Then the Banach space $\mathcal{H}^p(\cG,\tilde \cG)$ is defined to be
\begin{equation}
\label{Hpwdef}
\mathcal{H}^p(\cG,\tilde{\cG}) := \{ f \in \mathcal{H} : \quad f= \sum_{n \in {\mathbb{Z}}} \langle f, \tilde g_n \rangle g_n, \quad (\langle f,\tilde g_n\rangle)_{n \in {\mathbb{Z}}} \in \ell^p({\mathbb{Z}})\}
\end{equation}
with the norm $\|f\|_{\mathcal{H}^p}= \| (\langle f,\tilde g_n \rangle )_{n \in {\mathbb{Z}}}\|_{\ell^p}$. 
Since $\ell^p({\mathbb{Z}}) \subset \ell^2({\mathbb{Z}})$, $\mathcal{H}^p$ is a dense subspace of $\mathcal{H}$. If $2 < p < \infty$ then we define $\mathcal{H}^p$ to be the completion of the subspace $\mathcal{H}_0$ of all finite linear combinations in $\cG$ with respect to the norm $\|f\|_{\mathcal{H}^p}= \| (\langle f,\tilde g_n \rangle )_{n \in {\mathbb{Z}}}\|_{\ell^p}$. If $p=\infty$ then we define $\mathcal{H}_0^\infty$ as the completion of $\mathcal{H}_0$ with respect to the norm $\|f\|_{\mathcal{H}^\infty_0}= \| (\langle f,\tilde g_n \rangle )_{n \in {\mathbb{Z}}}\|_{C_0}$ and $\cH^\infty:=(\cH^1)'$ endowed with the norm $\|f\|_{\mathcal{H}^\infty}:=\|(\langle f , \tilde g_n \rangle )_{n \in {\mathbb{Z}}}\|_{\ell^\infty}$.\\

\begin{rem}
The definition of $\mathcal{H}^p(\cG,\tilde{\cG})$ does not depend on the particular $\cA^1_s$-self localized dual chosen, and any other $\cA^1_s$-self-localized frame $\cF$ which is localized to $\cG$ generates in fact the same spaces. In particular one has that $\cH^2(\cG,\tilde \cG)=\cH$.
\end{rem}

Then the following theorem \cite{FG} holds.

\begin{tm}
\label{banachfr}
Assume that $\cG$ is an $\cA^1_s$-self-localized frame for $\mathcal{H}$, for $s>0$. Then $\cG$ and its canonical dual frame $\tilde \cG$ are frame atomic decomposition for $\mathcal{H}^p(\cG,\tilde{\cG})$ with $C_{\tilde \cG}$ and $C_{\cG}$ as corresponding and respective coefficient operators.
\end{tm}

\section{Perturbation and localization of frames}

In the last years several results on {\it stability} and {\it perturbation} of frames and Riesz bases have been investigated, for example, by Casazza, Christensen, and Heil \cite{C2,cc1,CH,FS}, to name some of the most prominent authors.
\\

One of the first motivations and classical applications of the perturbation of frame theory is the study of ``non-uniform'' {\it coherent} frames generated by a strongly continuous (square-integrable) and irreducible representation $\pi$ of some locally compact group by $\mathcal{U}(L^2(\mathbb{R}^d))$ (see for example \cite{FG3,G0}), i.e., $\mathcal{F}:=\{\pi(x)g\}_{x \in H}$, where $H \subset G$. In most of the classical cases, namely Gabor and wavelet frames, there is usually a well structured canonical choice of the index subset $H$, maybe a discrete subgroup. 
The question is whether, given a frame $\mathcal{F}:=\{\pi(x)g\}_{x \in H}$, a {\it perturbation} of $H$, namely $\tilde H$, might preserve the frame property, i.e., whether $\mathcal{E}:=\{\pi(\tilde x)g\}_{\tilde x \in \tilde H}$ is again a frame for $L^2(\mathbb{R}^d)$.
\\

Here we want to emphasize that the perturbation of frame theory is indeed very relevant and important for numerical purposes. In fact, it is never possible to compute exactly a frame (of functions), in particular its canonical dual frame, and the numerical methods applied in this context perform approximations up to some prescribed tolerance, and  numerical rounding errors are anyway present.
Then it is clear that if such perturbations destroyed the frame property, then the use of such expansions for numerical application would be potentially incorrect.
Moreover, in many applications, for example in PDE numerical solution \cite{S,DFR1}, it is relevant that the frames used can be stable in a wider sense, i.e., any small perturbation should preserve not only the Hilbert space frame properties, but even the frame atomic decomposition one, especially in the characterization of those Banach spaces where it is expected that the solution is sitting.
\\

Thus, in order to model what happens in concrete numerical approximations, here we want to discuss the following problems: 
\begin{itemize}
\item[1)] Assume that   $\cF:=\{f_n\}_{n \in \mathbb{Z}}$ and $\cG:=\{g_x\}_{x \in \mathbb{Z}}$ are two $\cA_s^1$-intrinsically localized frames in $\cH$ for $s>0$, and that a system $\cE=\{e_n\}_{n\in \mathbb{Z}}$ in $\cH$ has the following property:
\begin{equation}
\label{unif}
\|e_n -f_n \|_{\cH^p(\cG,\tilde \cG)} \leq \varepsilon_n, \quad \text{ for } 1 \leq p \leq \infty,
\end{equation}
for a positive sequence $\{\varepsilon_n\}_{n \in \mathbb{Z}}$ of tolerance errors. 
Are there conditions implying that $\cE$ is again a frame for $\cH$? 
\item[2)] Moreover, if $\cF \sim_{\cA^1_s} \cG$ then it is known (see Remark in the previous section) that $\cH^p:=\cH^p(\cF,\tilde \cF)=\cH^p(\cG,\tilde \cG)$ and that $\cF$ and $\tilde \cF$ are frame atomic decomposition for $\cH^p$. Are there also  conditions implying that $\cE$ is again a frame atomic decomposition for $\cH^p$?
\end{itemize}

\begin{rem}
In \eqref{unif} we assume an approximation in $\cH^p$. The motivation is that the so-called non-linear or best $N$-term approximation has usually convergence in such spaces (for all $p_0<p\leq\infty$, for some $p_0<1$) and it is performed by adaptive algorithms, see for example \cite{DFR1}. Thus, for later use, we will limit our discussion assuming such norm approximation, even if the arguments here illustrated can be developed with approximations considered in more general Banach spaces, see \cite{FG}.
\end{rem}

Let us discuss the Problem 1).
One of the first results on perturbation of frames has been proposed by Christensen in \cite{C2} and it reads as follows:

\begin{tm}
Let $\cF$ be a frame for $\cH$ with bounds $A,B$ and let $\cE$ be a system in $\cH$. If
\begin{equation}
\label{ch1}
        R:=\sum_{n\in \mathbb{Z}} \|e_n-f_n\|_{\cH}^2 < A,
\end{equation}
then $\cE$ is again a frame with bounds $A(1 - \sqrt{\frac{R}{A}})^2, B(1 - \sqrt{\frac{R}{B}})^2$.
\end{tm}
In our setting this result can be re-formulated in the following way: If $\sum_{n \in \mathbb{Z}} \varepsilon_n^2 < A$ then  $\cE$ is again a frame. Unfortunately this means that the approximation tolerated for each individual element of the frame is not uniform but should decrease quickly for $n$ going to $\infty$. 
For frames ``uniformly distributed'' in the space and well structured, one would prefer not to {\it treat} differently one element of the frame with respect to others, but would accept a more ``democratic'' (uniform) approximation.  
For this reason we propose here a slightly weaker condition based on an estimation of the $\cH^\infty$ norm of $(e_n-f_n)$.

\begin{prop}
\label{for}
Let $\cF$ be a frame for $\cH$ with bounds $A,B$ and let $\cE$ be a system in $\cH$. If 
\begin{equation}
\label{fo1}
\|( \|e_n- f_n\|_{\cH^1} )_{n\in\mathbb{Z}} \|_{\ell^\infty} \leq \varepsilon \text{ and } \|( \|e_n- f_n\|_{\cH^\infty} )_{n\in\mathbb{Z}} \|_{\ell^1} \leq  \varepsilon
\end{equation}
for $0 < \varepsilon< \frac{1}{\sqrt{A^{-1}}}$, then $\cE$ is again a frame with bounds $ A\L 1- \sqrt{A^{-1}} \varepsilon \R^2$, $ B\L1 + \sqrt{B^{-1}} \varepsilon\R^2$.
\end{prop}
\begin{proof}
Consider the coefficient map $C_{\cE-\cF}:\cH^\infty_0 \rightarrow C_0$, given by $C_{\cE-\cF}(f) = (\langle f, e_n- f_n\rangle)_{n \in \mathbb{Z}}$. This map is in fact bounded
$$
\sup_{n \in \mathbb{Z}} | \langle f, e_n- f_n\rangle| \leq \|f\|_{\cH^\infty_0} \sup_{n \in \mathbb{Z}} \|e_n- f_n\|_{\cH^1} \leq  \varepsilon \|f\|_{\cH^\infty_0}.
$$
This implies that $C_{\cE-\cF}^*:\ell^1 \rightarrow \cH^1$, given by $C_{\cE-\cF}^*(\mathbf{c})=\sum_{n\in \mathbb{Z}} c_n (e_n-f_n)$, is also bounded by $ \varepsilon$. Moreover 
$$
\|C_{\cE-\cF}^*(\mathbf{c})\|_{\cH^\infty} \leq \| \mathbf{c}\|_{\ell^\infty} \sum_{n \in \mathbb{Z}} \|e_n-f_n\|_{\cH^\infty} \leq \varepsilon  \|\mathbf{c}\|_{\ell^\infty}.
$$
This means that $C_{\cE-\cF}^*$ is also bounded by $\varepsilon$ from $\ell^\infty$ to $\cH^\infty$. By complex interpolation, one has that $C_{\cE-\cF}^*$ is also bounded from $\ell^2$ to $\cH$ with bound $\varepsilon$. This implies that
$$
C_{\cE} :\ell^2 \rightarrow \cH, \quad C_{\cE}(f)=  (\langle f, e_n\rangle)_{n \in \mathbb{Z}}
$$ 
is a bounded operator with bound $(\sqrt{B}+ \varepsilon)$, and, in particular that $\cE$ is a Bessel sequence. Moreover, observe that the operator 
$$
T f = \sum_{n \in \mathbb{Z}} \langle f, \tilde f_n \rangle e_n,
$$
is bounded and 
$$
\|(I-T) f \| = \|  \sum_{n \in \mathbb{Z}} \langle  f, \tilde f_n \rangle (e_n-f_n) \| \leq  \varepsilon \|(\langle  f, \tilde f_n \rangle)_{n \in \mathbb{Z}} \|_{\ell^2} \leq \sqrt{A^{-1} } \varepsilon\|f \|.
$$
Since $ \sqrt{A^{-1}} \varepsilon < 1$ then $T$ is an invertible operator and
$$
f = T T^{-1} f = \sum_{n \in \mathbb{Z}} \langle T^{-1} f, \tilde f_n \rangle e_n
$$
This implies that
$$
\|f\|^4 = \langle f, f \rangle^2 \leq \L \sum_{n \in \mathbb{Z}} |\langle T^{-1} f, \tilde f_n \rangle| \langle f, e_n \rangle| \R^2 \leq A^{-1} \|T^{-1}\|^2 \|f\|^2 \sum_{n \in \mathbb{Z}}| \langle f, e_n \rangle|^2
$$
and 
$$
 A \|T^{-1}\|^{-2} \|f\|^2 \leq \sum_{n \in \mathbb{Z}}| \langle f, e_n \rangle|^2.
$$
Since $\|T^{-1}\| \leq \frac{1}{1- \sqrt{A^{-1}} \varepsilon}$ one has that $\cE$ is a again a frame with bounds $ A\L 1- \sqrt{A^{-1}} \varepsilon \R^2$,  $ B\L1 + \sqrt{B^{-1}} \varepsilon\R^2$.
\end{proof}

On one hand this criterion appears weaker then \eqref{ch1} since  $\|\cdot\|_{\cH^\infty} \leq \|\cdot \|_{\cH^2}$. On the other hand, the approximation required at the level of the $\|\cdot\|_{\cH^\infty}$ norm is still not uniform, but it must anyway decrease for $n$ going to $\infty$. \\

A very strong result has been achieved by  Casazza and Christensen  which overcomes this difficulty. A very general version of this result can be found in \cite{cc1}.
\begin{tm}
\label{cc}
Let $\cF$ be a frame for $\cH$ with bounds $A,B$. Let $\cE$ be a system in $\cH$, and assume that there exist constants $\lambda, \mu \geq 0$ such that $\lambda + \frac{\mu}{\sqrt{A}} <1$ and
\begin{equation}
\label{ch2}
\|\sum_{n \in F} c_n (e_n - f_n) \| \leq \lambda \|\sum_{n \in F} c_n f_n\| + \mu \L \sum_{n \in F} |c_n|^2\R^{1/2},
\end{equation}
for any finite sequence $(c_n)_{n \in F}$ of scalars and $F \subset \mathbb{Z}$, $\# F < \infty$. Then $\cE$ is a frame with bounds $A(1 - (\lambda + \frac{\mu}{\sqrt{A}}))^2, B(1+ \lambda + \frac{\mu}{\sqrt{B}})^2$.
\end{tm}

Criterion \eqref{ch1} implies \eqref{ch2}, and, in particular, the latter does not imply that the approximation should improve for $n$ going to $\infty$.
For general frames, to verify the condition \eqref{ch2} might be anyway difficult. 
We refer the reader to the works \cite{F1,F3,F2,F4} where perturbation techniques similar to \eqref{ch2} have been exploited in order to prove the existence of a large class of intermediate (wave-packet) frames between the more classical and well-known Gabor and wavelet frames.

We want to illustrate here some sufficient conditions, usually easier to check for localized frames systems, which ensure the application of Theorem \ref{cc} as a useful tool in several concrete cases. 

\begin{prop}
\label{for2}
Let $\cF$ be a frame for $\cH$ with bounds $A,B$ and let $\cE$ be a system in $\cH$  satisfying 
\begin{equation}
\label{fo2}
\|A(\cE-\cF,\tilde \cG)\|_{\cA^1}= \max\{       \sup_{n \in \mathbb{Z}} \sum_{x \in \mathbb{Z}} |\langle e_n- f_n, \tilde g_x \rangle |,  \sup_{x \in \mathbb{Z}} \sum_{n \in \mathbb{Z}} |\langle e_n- f_n, \tilde g_x \rangle | \} \leq \varepsilon,
\end{equation}
with $0< \varepsilon < (\sqrt{A^{-1}\|A(\cG,\cG)\|_{\cA^1}})^{-1}$.
Then  $\cE$ is a frame with bounds $A(1- \sqrt {A^{-1} \| A(\cG,\cG)\|_{\cA^1}} \varepsilon)^2$, $B(1 + \sqrt {B^{-1}\| A(\cG,\cG)\|_{\cA^1}} \varepsilon)^2$.
\end{prop}

\begin{proof}
It is sufficient to show that $C^*_{\cE-\cF}$ is bounded from $\ell^2$ to $\cH$ with bound small enough. 
\begin{eqnarray*} 
\|C^*_{\cE-\cF}(\mathbf{c})\|^2 &=& \|\sum_{n \in \mathbb{Z}} c_n (e_n -f_n) \|^2 \\
&=& \sum_{n,m \in \mathbb{Z}} c_n \overline{c_m} \langle  e_n -f_n ,  e_m -f_m \rangle \\
&=& \sum_{n,m \in \mathbb{Z}} c_n \overline{c_m} \sum_{x,y \in \mathbb{Z}} \langle  e_n -f_n ,\tilde g_x \rangle \overline{\langle  e_m -f_m ,\tilde g_y \rangle} \langle g_x, g_y \rangle.
\end{eqnarray*}
Since $\cG$ is assumed $\cA^1_s$-intrinsically localized one has $A(\cG,\cG) \in \cA^1_s$.
Moreover, since  \eqref{fo2} holds one also has that $A(\cE-\cF,\cG) \in \cA^1$ with norm $\| A(\cE-\cF,\cG)\|_{\cA^1} \leq  \varepsilon$. 
By the property (A0) of the Schur algebra one immediately has
$$
\|C^*_{\cE-\cF}(\mathbf{c})\|^2 \leq \|\mathbf{c}\|_{\ell^2}^2 \| A(\cG,\cG)\|_{\cA^1}  \varepsilon^2.
$$
Therefore, since $\sqrt{A^{-1} \| A(\cG,\cG)\|_{\cA^1}}  \varepsilon <1$, then one concludes as in the proof of Proposition \ref{for}.
\end{proof}

Observe that if $\cF \sim_{\cA^1_s} \cG$ then the assumptions of Proposition \ref{for2} in particular imply that $\cE \sim_{\cA^1} \cG$. In fact, the localization principle has been thought and invented as a measure of similarity or equivalence of frames. Thus, it is not surprising that a natural localization condition \eqref{fo2} can be used as a measure of perturbation for a frame.

Let us summarize the result implications of this section in the following theorem.

\begin{tm}
Let $\cF$ be a frame for $\cH$  and let $\cE$ be a system in $\cH$  satisfying \eqref{unif} with $\sup_{n \in \mathbb{Z}} \varepsilon_n \leq \varepsilon$.
Then the following statements are ordered by implications
\begin{itemize}
\item[i)]  $\sum_{n\in \mathbb{Z}}  \|e_n-f_n\|_{\cH^\infty}  \leq \varepsilon$;
\item[ii)] $\sup_{x \in \mathbb{Z}} \sum_{n \in \mathbb{Z}} |\langle e_n- f_n, \tilde g_x \rangle | \leq \varepsilon$;
\item[iii)] $ \|\sum_{n \in \mathbb{Z}} c_n (e_n -f_n) \| \leq  \sqrt{\| A(\cG,\cG)\|_{\cA^1}} \varepsilon \L \sum_{n \in F} |c_n|^2\R^{1/2}$.
\end{itemize}
In particular if $\varepsilon>0$ is small enough and one of i)-iii) is valid then $\cE$ is again a frame for $\cH$.
\end{tm}

We want to conclude this paper with a stability result of $\cE$ as a frame atomic decomposition for $\cH^p$ for all $p \in [1,\infty]$, and as an answer to Problem 2).
\begin{tm}\label{for3}
Let $\cF \sim_{\cA^1_s} \{\cF,\cG\}$. This implies that  $\cF$ is an atomic decomposition for $\cH^p$. Let us denote the atomic decomposition bounds $A_p,B_p>0$ for all $p \in [1, \infty]$, and assume $B:=\sup_{p \in [1,\infty]} B_p<\infty$. 
Let $\cE$ be a system in $\cH$ satisfying
\begin{equation}
\label{fo3}
\|A(\cE-\cF,\tilde \cG)\|_{\cA^1}= \max\{        \sup_{n \in \mathbb{Z}} \sum_{x \in \mathbb{Z}} |\langle e_n- f_n, \tilde g_x \rangle |, \sup_{x \in \mathbb{Z}} \sum_{n \in \mathbb{Z}} |\langle e_n- f_n, \tilde g_x \rangle | \}\leq  \varepsilon,
\end{equation}
for $0<\varepsilon<B^{-1}$.
Then  $\cE\sim_{\cA^1} \{\cE, \cF,\cG\}$ is a frame atomic decomposition for $\cH^p$  with bounds $A_p(1+(\varepsilon  B_p))^{-1}$, $B_p(1-(\varepsilon  B_p))^{-1}$ for all $p \in [1,\infty]$.
 
\end{tm}

\begin{proof}
By similar arguments as in the proof of Proposition \ref{for2} one can show that
$$ \|\sum_{n \in \mathbb{Z}} c_n (e_n -f_n) \|_{\cH^p} \leq   \varepsilon \L \sum_{n \in \mathbb{Z}} |c_n|^p\R^{1/p}.
$$
By an application of \cite[Theorem 2.3]{CH} one has that $\cE$ is a frame atomic decomposition for $\cH^p$  with bounds $A_p(1+(\varepsilon  B_p))^{-1}$, $B_p(1-(\varepsilon  B_p))^{-1}$ for all $p \in [1,\infty]$. Since $\cE \sim_{\cA^1} \cG$ and $\tilde \cG \sim_{\cA^1} \cF$ then $\cE \sim_{\cA^1} \cF$.
It remains to show that $\cE\sim_{\cA^1} \cE$. Note that
$$
|A(\cE,\cE)| \leq |A(\cE,\tilde \cG)||A(\cG,\cE)|.
$$
All the matrices on the right-hand side are in $\cA^1$ and then by (A3) one has that $A(\cE,\cE) \in \cA^1$. 
\end{proof}

This result is quite interesting since, for $s=0$ Schur type localization, it is not known whether there exists dual frames  $\tilde \cE$ with such localization to define the corresponding spaces $\cH^p(\cE,\tilde \cE)$. Therefore, Theorem \ref{for3} gives,  by a perturbation argument, an alternative criterion to Theorem \ref{self} in order to show that a $\cA^1$-intrinsically localized frame $\cE$ in fact can extend to a frame atomic decomposition for a class of non-trivial Banach spaces $\cH^p$. \\

\noindent {\bf ACKNOWLEDGEMENT:} The author would like to thank Hans G. Feichtinger and Karlheinz Gr\"ochenig for the valuable discussions and insights, and the hospitality of NuHAG (the Numerical Harmonic Analysis Group), Dept. Maths, University of Vienna, during the preparation of this work.

\bibliography{banach} 

\noindent Massimo Fornasier\\
Universit\`a ``La Sapienza'' in Roma\\
Dipartimento di Metodi e Modelli Matematici per le Scienze Applicate\\
Via Antonio Scarpa, 16/B\\
00161 Roma\\
Italy\\
email: {\tt mfornasi@math.unipd.it}\\
WWW: {\tt http://www.math.unipd.it/$\sim$mfornasi}

\end{document}